\documentclass[fleqn]{mat01}
\usepackage{times,mathtimy,amssymb,latexsym}%,krepsf}
\begin{document}

\setcounter{page}{453}
\firstpage{453}

\def\d{\hbox{d}}
\def\e{\hbox{e}}

\newtheorem{theore}{Theorem}
\newtheorem{theor}[theore]{\bf Theorem}
\newtheorem{definit}[theore]{\rm DEFINITION}
\newtheorem{lem}[theore]{Lemma}
\newtheorem{rem}[theore]{Remark}
\newtheorem{propo}[theore]{\rm PROPOSITION}
\newtheorem{coro}[theore]{\rm COROLLARY}

\title{Topologically left invariant means on semigroup algebras}

\markboth{Ali Ghaffari}{Topologically left invariant means on
semigroup algebras}

\author{ALI GHAFFARI}

\address{Department of Mathematics, Semnan University, Semnan, Iran\\
\noindent E-mail: ghaffari1380@yahoo.com}

\volume{115}

\mon{November}

\parts{4}

\pubyear{2005}

\Date{MS received 22 February 2005; revised 17 July 2005}

\begin{abstract}
Let $M(S)$ be the Banach algebra of all bounded regular Borel
measures on a locally compact Hausdorff semitopological semigroup
$S$ with variation norm and convolution as multiplication. We
obtain necessary and sufficient conditions for $M(S)^{*}$ to have
a topologically left invariant mean.
\end{abstract}


\maketitle

\setcounter{equation}{0}

\section{Introduction}

Let $S$ be a locally compact Hausdorff semitopological semigroup
with convolution measure algebra $M(S)$ and probability measures
$M_{0}(S)$. We know that $M(S)$ is a Banach algebra with total
variation norm and convolution. The first Arens multiplication on
$M(S)^{*}$ is defined in three steps as follows.

For $\mu, \nu$ in $M(S), f$ in $M(S)^{*}$ and $F, G$ in
$M(S)^{**}$, the elements $f\mu, Ff$ of $M(S)^{*}$ and $GF$ of
$M(S)^{**}$ are defined by
\begin{equation*}
\langle f\mu, \nu \rangle = \langle f, \mu * \nu \rangle, \quad
\langle Ff, \mu \rangle = \langle F, f\mu \rangle, \quad \langle
GF, f \rangle = \langle G, Ff \rangle.
\end{equation*}
Denote by 1 the element in $M(S)^{*}$ such that $\langle 1, \mu
\rangle = \mu (S), \mu \in M(S)$. A linear functional $M \in
M(S)^{**}$ is called a {\it mean} if $\langle M, f \rangle\geq 0$
whenever $f\geq 0$ and $\langle M, 1\rangle = 1$. Each probability
measure $\mu \in M_{0}(S)$ is a mean. An application by the
Hahn--Banach theorem shows that $M_{0}(S)$ is weak$^{*}$ dense in
the set of means on $M(S)^{*}$. A mean $M$ is {\it topological
left invariant} if $\langle M, f\mu \rangle = \langle M, f
\rangle$ for any $\mu \in M_{0}(S)$ and $f\in M(S)^{*}$. We shall
follow Ghaffari \cite{7} and Wong \cite{14,15} for definitions and
terminologies not explained here. We know that topologically left
invariant mean on $M(S)^{*}$ have been studied by Riazi and Wong
in \cite{11} and by Wong in \cite{14,15}. They also went further
and for several subspaces $X$ of $M(S)^{*}$, have obtained a
number of interesting and nice\break results.

The existence of topologically left invariant means and left
invariant means for groups was investigated widely by Paterson
\cite{9} and Pier \cite{10}. Other important studies on amenable
semigroups are those of Argabright \cite{1}, Day \cite{4}, Lau
\cite{6}, and Mitchell \cite{8}. For further studies and
complementary historical comments see \cite{3,9,10}.

Let $M_{0}(S)$ have a measure $\nu$ such that the map $s \mapsto
\delta_{s} * \nu$ from $S$ into $M(S)$ is continuous. The author
recently proved that the following conditions are equivalent:\pagebreak
\begin{enumerate}
\leftskip .15pc
\renewcommand\labelenumi{\rm (\alph{enumi})}

\item $M(S)^{*}$ has a topologically left invariant mean;
\item there is a net $(\mu_{\alpha})$ in $M_{0}(S)$ such that for
every compact subset $K$ of $S$, $\|\mu * \mu_{\alpha} -
\mu_{\alpha}\|$ $\rightarrow 0$ uniformly over all $\mu$ in
$M_{0}(S)$ which are supported in $K$.\vspace{-.5pc}
\end{enumerate}

In this paper, we obtain a necessary and sufficient condition for
$M(S)^{*}$ to have a topologically left invariant mean.

\section{Main results}

Throughout the paper, $S$ is a locally compact Hausdorff
semitopological semigroup. We say that $S$ is {\it semifoundation}
if there is a measure $\nu \in M_{0}(S)$ such that the map
$x\mapsto \delta_{x} * \nu$ from $S$ into $M(S)$ is continuous. It
is clear that every foundation semigroup is also a semifoundation
semigroup (for more on foundation semigroups, the reader is
referred to \cite{2} and \cite{5}). We recall that a mean $M$ is
{\it left invariant} if $\langle M, f\delta_{x} \rangle = \langle
M, f \rangle$ for any $x\in S$ and $f\in M(S)^{*}$. Obviously, a
topologically left invariant mean on $M(S)^{*}$ is also a left
invariant mean on $M(S)^{*}$.

\begin{propo}$\left.\right.$\vspace{.5pc}

\noindent Let $S$ be a semifoundation semigroup. Choose $\nu \in
M_{0}(S)$ such that the map $x\mapsto \delta_{x} * \nu$ from $S$
into $M(S)$ is continuous. If $M\in M(S)^{**}$ is a left invariant
mean on $M(S)^{*}$ and $\langle M, f\nu \rangle = \langle M, f
\rangle$ for each $f\in M(S)^{*}${\rm ,} then $M$ is a
topologically left invariant mean on $M(S)^{*}$.
\end{propo}

Note that this is Proposition~22.2 of \cite{10}, which was proved
for groups. However, our proof is completely different.

\begin{proof}
Suppose that $f\in M(S)^{*}$ and $\mu \in M_{0}(S)$. For every
$x\in S$, we can write
\begin{equation*}
\langle M, f\delta_{x} * \nu \rangle = \langle M,
f\delta_{x}\rangle = \langle M, f \rangle.
\end{equation*}
It follows that
\begin{equation}
\int \langle M f, \delta_{x} *\nu \rangle \d\mu (x) = \langle M, f
\rangle.
\end{equation}
Since $x\mapsto\delta_{x} * \nu$ is continuous, by Theorem~3.27 in
\cite{12}, it is easy to see that
\begin{equation}
\int \langle Mf, \delta_{x} * \nu \rangle\d \mu (x) = \langle Mf,
\mu * \nu \rangle = \langle M, f\mu * \nu \rangle.
\end{equation}
Hence, using (1) and (2), $\langle M, f\mu * \nu \rangle = \langle
M, f \rangle$. By hypothesis,
\begin{equation*}
\langle M, f\mu * \nu\rangle = \langle M, (f\mu) \nu \rangle =
\langle M, f\mu \rangle.
\end{equation*}
Consequently $\langle M, f \rangle = \langle M, f\mu \rangle$,
i.e., $M$ is a topologically left invariant mean on\break
$M(S)^{*}$. \hfill $\Box$
\end{proof}

Our next result gives an important property that characterizes
topologically left amenability of $M(S)^{*}$.

\setcounter{theore}{0}

\begin{theor}[\!]
Let $S$ be semifoundation semigroup with identity. Then the
following statements are equivalent{\rm :}

\begin{enumerate}
\leftskip .15pc
\renewcommand\labelenumi{{\rm (\arabic{enumi})}}
\item $M(S)^{*}$ has a topologically left invariant mean{\rm ;}
\item for all $n\in {\mathbb N}$ and $\mu_{1},\ldots, \mu_{n}\in
M(S),$
\begin{align*}
\hskip -1.25pc &\inf \{\sup \{\|\mu_{i} * \mu\|; 1 \leq i \leq n
\}, \mu \in M_{0}(S)\} \leq \sup \{|\mu_{i}(S)|; 1 \leq i \leq
n\}.
\end{align*}
\end{enumerate}
\end{theor}

\begin{proof}
Let $M(S)^{*}$ have a topologically left invariant mean. Let
$\mu_{1},\ldots, \mu_{n}\in$ $M(S)$, $\epsilon > 0$ and put
$\delta = \epsilon (2 + \sup \{\|\mu_{i}\|; 1 \leq i \leq
n\})^{-1}$. There exists a compact subset $K$ in $S$ such that
$|\mu_{i}|(S \backslash K) < \delta$ whenever $i = 1, \ldots, n$.
By Theorem~2.2 in \cite{7}, there exists a measure $\mu$ in
$M_{0}(S)$ such that the map $x \mapsto \delta_{x} * \mu$ from $S$
into $M(S)$ is continuous and $\|\delta_{x} * \mu - \mu \| <
\delta$ for any $x \in K$. Thus, for every $i = 1, \ldots, n$ and
$f\in M(S)^{*}$, by Theorem~3.27 in \cite{12}, we can write
\begin{align*}
|\langle f, \mu_{i} * \mu - \mu_{i}(S) \mu\rangle| &= \left|\int
\langle f, \delta_{x} * \mu \rangle \d\mu_{i}(x) - \mu_{i}
(S)\langle f, \mu\rangle \right|\\[.4pc]
&= \left|\int (\langle f, \delta_{x} * \mu \rangle - \langle f,
\mu \rangle)\d \mu_{i}(x) \right|\\[.4pc]
&= \left|\int \langle f, \delta_{x} * \mu - \mu \rangle \d
\mu_{i}(x)\right|\\[.4pc]
&= \left|\int_{S\backslash K}\langle f, \delta_{x} * \mu - \mu
\rangle \d\mu_{i}(x) \right.\\[.4pc]
&\quad\, \left. + \int_{K}\langle f, \delta_{x} * \mu - \mu
\rangle \d\mu_{i}(x) \right|\\[.4pc]
&\leq 2|\mu_{i}| (S\backslash K)\|f\| + \delta \|f\||\mu_{i}|(K)\\[.4pc]
&\leq 2\delta\|f\| + \delta \|f\|\|\mu_{i}\|\\[.35pc]
&= \delta\|f\| (2 + \|\mu_{i}\|) \leq \epsilon \|f\|.
\end{align*}
It follows that $\|\mu_{i} * \mu - \mu_{i}(S)\mu\| < \epsilon$
whenever $i = 1, \ldots, n$. Consequently
\begin{equation*}
\sup \{\|\mu_{i} * \mu\|; 1 \leq i \leq n\} \leq \sup
\{|\mu_{i}(S)|; 1\leq i \leq n \} + \epsilon.
\end{equation*}
Therefore
\begin{equation*}
\inf \{\sup \{\|\mu_{i} * \mu\|; 1 \leq i \leq n\}; \mu \in
M_{0}(S)\}\leq \sup\{|\mu_{i}(S)|; 1 \leq i \leq n\}.
\end{equation*}

Conversely let $\mu_{1},\ldots, \mu_{n}\in M_{0}(S)$ and $\epsilon
> 0$. For any $i = 1,\ldots, n$, consider $\nu_{i} = \mu_{i} -
\delta_{e}$. We have $\nu_{i}(S) = 0$ whenever $i = 1,\ldots, n$.
By assumption,
\begin{equation*}
\inf\{\sup\{\|\nu_{i} * \mu\|; 1 \leq i \leq n\}; \mu \in
M_{0}(S)\} = 0.
\end{equation*}
Thus there exists $\mu \in M_{0}(S)$ such that
\begin{equation*}
\sup\{\|\nu_{i} * \mu\|; 1 \leq i \leq n\} < \epsilon,
\end{equation*}
i.e., for every $i = 1,\ldots, n, \|\mu_{i} * \mu - \mu\| <
\epsilon$. By Theorem~2.2 in \cite{7}, $M(S)^{*}$ has a
topologically left invariant mean. \hfill $\Box$
\end{proof}

Let $V$ be a locally convex Hausdorff topological vector space and
let $Z$ be a compact convex subset of $V$. The pair $(M_{0}(S), Z)$ is
called a {\it semiflow}, if;

\begin{enumerate}
\leftskip .15pc
\renewcommand\labelenumi{\rm (\arabic{enumi})}
\item There exists a map $\rho\hbox{\rm :}\ M(S) \times V \rightarrow V$ such
that for every $z\in Z$, the map $\rho (-, z)\hbox{\rm :}\
M(S)\rightarrow V$ is continuous and linear $(M(S)$ has the
topology $\sigma (M(S)$, $M(S)^{*}))$;
\item $\rho(M_{0}(S), Z) \subseteq Z$;
\item For any $\mu, \nu\in M (S)$ and $z\in Z, \rho(\mu, \rho(\nu, z)) =
\rho(\mu * \nu, z)$.
\end{enumerate}
We remind the reader of our notation conventions:
\begin{equation*}
\mu z = \rho(\mu, z), \quad \mu\in M(S), z \in Z.
\end{equation*}

\begin{theor}[\!]
Let $S$ be a semitopological semigroup. The following statements
are equivalent{\rm :}

\begin{enumerate}
\leftskip .15pc
\renewcommand\labelenumi{\rm (\arabic{enumi})}
\item $M(S)^{*}$ has a topologically left invariant mean{\rm ;}
\item for every $f\in M(S)^{*}${\rm ,} there exists a mean $M$ such that
$\langle M, f\mu \rangle = \langle M, f\nu \rangle$ for any $\mu,
\nu$ in $M_{0}(S)${\rm ;}
\item \looseness 1 for any semiflow $(M_{0}(S), Z)${\rm ,} there is some $z\in Z$
such that $\mu z = z$ for all $\mu \in M_{0}(S)$.
\end{enumerate}
\end{theor}

\begin{proof}
(1) implies (2) is easy.

Now, assume that (2) holds. We will show that $M(S)^{*}$ has a
topologically left invariant mean. To each $f\in M (S)^{*}$, we
associate the non-void subset
\begin{equation*}
\Omega_{f} = \{M \in \Omega; \langle M, f\mu\rangle = \langle M,
f\nu\rangle \quad \hbox{for all} \ \mu, \nu \in M_{0}(S)\},
\end{equation*}
($\Omega$ is the convex set of all means on $M(S)^{*}$.) The sets
$\Omega_{f}$ are obviously weak$^{*}$ compact. We shall show that
the family $\{\Omega_{f}; f\in M(S)^{*}\}$ has the finite
intersection property. Since $\Omega$ is weak$^{*}$ compact, it
will follow that
\begin{equation*}
\bigcap \{\Omega_{f}; f \in M(S)^{*}\} \neq \emptyset;
\end{equation*}
and if $M$ is any member of this intersection, then $M^{2}$ is a
topologically left invariant mean on $M(S)^{*}$.

We proceed by induction. By hypothesis, $\Omega_{f} \neq
\emptyset$ for each $f \in M(S)^{*}$. Let $n\in \mathbb{N},
f_{1},\ldots, f_{n}\in M(S)^{*}$ and assume that $\cap_{i = 1 }^{n
- 1}\Omega_{f_{i}}\neq \emptyset$. If $M_{1}$ is a member of this
intersection and if $M_{2}\in \Omega_{M_{1}f_{n}}$, then for every
$\mu, \nu$ in $M_{0}(S)$ we have
\begin{equation*}
\langle M_{2}M_{1},f_{n}\mu \rangle = \langle M_{2},
M_{1}f_{n}\mu\rangle = \langle M_{2}, M_{1} f_{n}\nu \rangle =
\langle M_{2}M_{1}, f_{n}\nu \rangle
\end{equation*}
and, for $i = 1, \ldots, n - 1$,
\begin{align*}
\langle M_{2}M_{1}, f_{i}\mu \rangle &= \langle M_{2},
M_{1}f_{i}\mu \rangle = \lim_{\alpha}\langle\mu_{\alpha},
M_{1}f_{i}\mu \rangle\\[.35pc]
&= \lim_{\alpha}\langle M_{1}f_{i}\mu, \mu_{\alpha}\rangle =
\lim_{\alpha}\langle M_{1}, (f_{i}\mu)\mu_{\alpha} \rangle\\[.35pc]
&= \lim_{\alpha}\langle M_{1},f_{i}\mu * \mu_{\alpha}\rangle =
\lim_{\alpha}\langle M_{1}, f_{i}\nu * \mu_{\alpha} \rangle\\[.35pc]
&= \lim_{\alpha}\langle M_{1},(f_{i}\nu) \mu_{\alpha}\rangle =
\lim_{\alpha}\langle \mu_{\alpha}, M_{1}f_{i}\nu \rangle\\[.35pc]
&= \langle M_{2}M_{1}, f_{i}\nu \rangle.
\end{align*}
(Recall that $M_{0}(S)$ is weak$^{*}$ dense in $\Omega$, and so
there is a net $\{\mu_{\alpha}\}$ in $M_{0}(S)$ such that
$\mu_{\alpha} \rightarrow M_{2}$ in the weak$^{*}$ topology.)
Hence $M_{2}M_{1}\in \cap_{i = 1}^{n}\Omega_{f_{i}}$. Thus
$\{\Omega_{f}; f\in M(S)^{*}\}$ has the finite intersection
property, as required. So (1) is equivalent to (2).

To prove that (1) and (3) are equivalent, let $(M_{0}(S), Z)$ be a
semiflow on a compact convex subset $Z$ of a locally convex
Hasudorff topological vector space $V$. If $f\in V^{*}$ and $z \in
Z$, we consider the mapping $f^{z}\hbox{\rm :}\ M(S)\rightarrow
\mathbb{C}$ given by $\langle f^{z}, \mu \rangle = \langle f, \mu
z \rangle$. It is easy to see that $f^{z}\in M(S)^{*}$. Let
$\Omega$ be the convex set of all means on $M(S)^{*}$. For $M \in
\Omega$, we can define $T(M)\hbox{\rm :}\ V^{*} \rightarrow
\mathbb{C}$ given by $\langle T(M), g\rangle = \langle M, g^{z}
\rangle$ $(g \in V^{*})$. One easily notes that $T(M)$ is linear.
Now, we embed $V$ into the algebraic dual $V^{*\prime}$ of $V^{*}$
with the topology $\sigma (V^{*\prime}, V^{*})$. Since $Z$ is
compact in $V$, it is closed in $V^{*\prime}$. On the other hand,
for every $h\in V^{*}$ and $\mu\in M_{0}(S)$, we have
\begin{equation*}
\langle T(\mu), h \rangle = \langle \mu, h^{z} \rangle = \langle
h, \mu z \rangle = \langle \mu z, h\rangle.
\end{equation*}
It follows that the $M_{0}(S)$-invariance of $Z$ implies that
$T(\mu) \in Z$. Since $M_{0}(S)$ is weak$^{*}$-dense in $\Omega$
and $Z$ is closed in $V^{*\prime}$, we conclude that $T(M)\in Z$
for every $M\in \Omega$. If $\mu \in M_{0}(S)$, we consider
$\lambda_{\mu}\hbox{\rm :}\ Z$ $\rightarrow Z$ by
$\lambda_{\mu}(z) = \mu z(z\in Z)$. Now let $M$ be a
topologically left invariant mean on $M(S)^{*}$. For every $h\in
V^{*}$ and $\mu \in M_{0}(S)$, we have
\begin{align*}
\langle \mu T(M), h \rangle &= \langle T(M), h \circ \lambda_{\mu}
\rangle = \langle M, (h \circ \lambda_{\mu})^{z}\rangle\\[.35pc]
&= \langle M, h^{z}\mu \rangle = \langle M, h^{z} \rangle\\[.35pc]
&= \langle T (M), h \rangle.
\end{align*}
So $\mu T(M) = T(M)$ for every $\mu \in M_{0}(S)$, i.e., $T(M)$ is
a fixed point under the action of $M_{0}(S)$.

To prove the converse, we know that the set $\Omega$ is convex and
weak$^{*}$-compact in $M(S)^{**}$. We define the semiflow
$(M_{0}(S), \Omega)$ by putting $\rho(\mu, F) = \mu F$ for $\mu
\in M(S)$ and $F\in M(S)^{**}$. By hypothesis, there exists $M\in
\Omega$ that is fixed under the action of $M_{0}(S)$, that is $\mu
M = M$ for every $\mu \in M_{0}(S)$. It follows that $M$ is a
topologically left invariant mean on $M(S)^{*}$. This completes
our proof. \hfill $\Box$
\end{proof}

A {\it right action} of $M(S)$ on $M(S)^{*}$ is a map $T\hbox{\rm
:}\ M(S) \times M(S)^{*} \rightarrow M(S)^{*}$ (denoted by $(\mu,
f) \mapsto T_{\mu}(f), \mu \in M(S)$ and $f\in M(S)^{*})$ such that

(1) $(\mu, f) \mapsto T_{\mu}(f)$ is bilinear and $T_{\mu
* \nu} = T_{\nu} \circ T_{\mu}$ for any $\mu, \nu \in M(S)$,

(2)
$T_{\mu}\hbox{\rm :}\ M(S)^{*}\rightarrow M(S)^{*}$ is a positive
linear operator and $T_{\mu}(1) = 1$ for any $\mu \in M_{0}(S)$.

Let $X$ be a linear subspace of $M(S)^{*}$ with $1 \in X$. We say
that $M \in X^{*}$ is a {\it mean} on $X$ if $\langle M, f
\rangle\geq 0$ if $f\geq 0$ and $\langle M, 1 \rangle = 1$. A mean
$M$ is $M_{0}(S)$-invariant under the right action $T$ if $\langle
M, T_{\mu}(f)\rangle = \langle M, f \rangle$ for any $\mu \in
M_{0}(S)$ and $f\in X$. We say that $X$ is $M_{0}(S)$-invariant
under the right action $T$ if $T_{\mu}(X)\subseteq X$ for any $\mu
\in M_{0}(S)$.

\begin{theor}[\!]
Let $S$ be a semitopological semigroup. The following satatements
are equivalent{\rm :}

\begin{enumerate}
\leftskip .15pc
\renewcommand\labelenumi{{\rm (\arabic{enumi})}}
\item $M(S)^{*}$ has a topologically left invariant mean{\rm ;}
\item for any separately continuous right action $T\hbox{\rm :}\ M(S) \times
M(S)^{*} \rightarrow M(S)^{*}$ of $M(S)$ on $M(S)^{*}$ {\rm (}$M(S)$ has
the topology $\sigma (M(S), M(S)^{*})$ and $M(S)^{*}$ has the weak
topology{\rm )} and any $M_{0}(S)$-invariant subspace $X$ of $M(S)^{*}$
containing $1,$ any $M_{0}(S)$-invariant mean $M$ on $X$ can be
extended to a $M_{0}(S)$-invariant mean $\cal{M}$ on $M(S)^{*}$.
\end{enumerate}
\end{theor}

\begin{proof}
Let $M(S)^{*}$ have a topologically left invariant mean, and let
\begin{equation*}
T\hbox{\rm :}\ M(S) \times M(S)^{*} \rightarrow M(S)^{*}
\end{equation*}
be a separately continuous right action of $M(S)$ on $M(S)^{*}$
and $M$ be a mean on $M_{0}(S)$-invariant subspace $X$ of
$M(S)^{*}$. Let
\begin{equation*}
Z = \{\cal{M}\in M(S)^{**}; \cal{M} \ \hbox{is a mean on} \
M(S)^{**} \ \hbox{and extends} \ M \}.
\end{equation*}
By the Hahn--Banach theorem, $Z\neq \emptyset$. It is easy to see
that $Z$ is a weak$^{*}$ closed convex subset of the unit ball in
$M(S)^{**}$, and is therefore weak$^{*}$ compact. Define
$\rho\hbox{\rm :}\ M(S)\times M(S)^{**} \rightarrow M(S)^{**}$ by
$\rho(\mu, F) = T_{\mu}^{*} (F), \mu \in M(S)$, $F \in M(S)^{**}$.
Notice that, since $T\hbox{\rm :}\ M(S) \times M(S)^{*}
\rightarrow M(S)^{*}$ is a separately continuous right action of
$M(S)$ on $M(S)^{*}$, it is clear that
\begin{equation*}
\rho (-, F)\hbox{\rm :}\ M(S) \rightarrow M(S)^{**}
\end{equation*}
is continuous for any $F\in M(S)^{**} (M(S)$ has the topology
$\sigma (M(S), M(S)^{*})$ and $M(S)^{**}$ has the topology $\sigma
(M(S)^{**}, M(S)^{*}))$. On the other hand, it is clear that each
$\rho(-, F)\hbox{\rm :}\ M(S) \rightarrow M(S)^{**}$ is linear
since $T\hbox{\rm :}\ M(S) \times M(S)^{*} \rightarrow M(S)^{*}$
is bilinear. Let $\cal{M}\in Z$ and $\mu \in M_{0}(S)$. Since
$T_{\mu}\hbox{\rm :}\ M(S)^{*} \rightarrow M(S)^{*}$ is positive
linear and $T_{\mu}(1) = 1$, so $T_{\mu}^{*}(\cal{M})$ is a mean
on $M(S)^{*}$. Now, let $f\in X$ we have
\begin{equation*}
\langle T_{\mu}^{*} (\cal{M}), f \rangle = \langle \cal{M},
T_{\mu} (f) \rangle = \langle M, T_{\mu}(f) \rangle = \langle M, f
\rangle.
\end{equation*}
This shows that $\rho(\mu, \cal{M}) = T_{\mu}^{*}(\cal{M})\in Z$,
i.e., $\rho(M_{0}(S), Z)\subseteq Z$. Let $\mu, \nu \in M(S)$ and
$\cal{M}\in Z$. Since $T\hbox{\rm :}\ M(S) \times M(S)^{*}
\rightarrow M(S)^{*}$ is an anti-homomorphism of $M(S)$ into the
algebra of linear operators in $M(S)^{*}$, therefore
\begin{align*}
\langle \rho(\mu, \rho (\nu, \cal{M})), f \rangle &= \langle
T_{\mu}^{*}(\rho (\nu, \cal{M})), f \rangle = \langle
T_{\mu}^{*}(T_{\nu}^{*}(\cal{M})), f \rangle\\[.3pc]
&= \langle T_{\nu}^{*}(\cal{M}), T_{u}(f)\rangle =
\langle\cal{M}, T_{\nu} (T_{\mu}(f))\rangle\\[.3pc]
&= \langle \cal{M}, T_{\mu *\nu}(f)\rangle = \langle T_{\mu * \nu}^{*}
(\cal{M}), f \rangle\\[.3pc]
&= \langle \rho (\mu * \nu, \cal{M}), f \rangle
\end{align*}
for any $f\in M(S)^{*}$. This shows that $\rho(\mu, \rho(\nu,
\cal{M})) = \rho(\mu * \nu, \cal{M})$ for any $\mu, \nu$ in $M(S)$
and $\cal{M}\in Z$. As we saw above, the pair $(M_{0}(S), Z)$ is a
semiflow. By Theorem~2, there is some $\cal{M} \in Z$ such that
$T_{\mu}^{*} (\cal{M}) = \rho(\mu, \cal{M}) = \cal{M}$ for each
$\mu \in M_{0}(S)$. $\cal{M}$ is then the required extension of
$M$.

Conversely, we define a right action $T\hbox{\rm :}\ M(S) \times
M(S)^{*} \rightarrow M(S)^{*}$ by putting $T_{\mu}(f) = f\mu$ for
$\mu \in M(S)$ and $f\in M(S)^{*}$. We claim that it is separately
continuous. If $\mu_{\alpha} \rightarrow \mu$ in the $\sigma
(M(S), M(S)^{*})$, then for any $F\in M(S)^{**}$, we have
\begin{align*}
\langle F, T_{\mu_{\alpha}}(f)\rangle &= \langle F,
f {\mu_{\alpha}} \rangle = \langle Ff, \mu_{\alpha} \rangle
\rightarrow \langle Ff, \mu \rangle\\[.3pc]
&= \langle F, f \mu \rangle = \langle F, T_{\mu}(f)\rangle.
\end{align*}
On the other hand, it is easy to see that every $T_{\mu}\hbox{\rm
:}\ M(S)^{*} \rightarrow M(S)^{*}$ is continuous ($M(S)^{*}$ has
the weak topology). Now choose $X$ to be the constants and define
$\langle M, \alpha \cdot 1 \rangle = \alpha$, for any $\alpha
\cdot 1 \in  X$. Then $M$ is a mean on $X$ satisfying $\langle M,
T_{\mu}(f) \rangle = \langle M, f\rangle$ for any $\mu\in
M_{0}(S)$ and $f\in X$. Any invariant extension $\cal{M}$ of $M$
to $M(S)^{*}$ is necessarily a topologically left invariant mean
on $M(S)^{*}$. This completes our proof. \hfill $\Box$
\end{proof}

The above characterization of topologically left invariant mean on
$M(S)^{*}$ is an analogue of Silverman's invariant extension
property in \cite{13}.

%\section*{Acknowledgements}
%The author is indebted to the University of Semnan for their
%support.

\end{document}